\documentclass[12pts]{amsart}

\usepackage{amscd,latexsym,amsthm,amsfonts,amssymb,amsmath,amsxtra}
\usepackage[mathscr]{eucal}
\newcommand{\BA}{{\mathbb {A}}}

\newcommand{\BC}{{\mathbb {C}}}

\newcommand{\GL}{{\mathrm{GL}}}

\newcommand{\Sp}{{\mathrm{Sp}}}

\newtheorem{thm}{Theorem}[section]
\newtheorem{cor}[thm]{Corollary}
\newtheorem{lem}[thm]{Lemma}
\newtheorem{prop}[thm]{Proposition}

\begin{document}
	\renewcommand{\theequation}{\arabic{equation}}
	\numberwithin{equation}{section}

	\title{Siegel-Weil type identities for general linear groups related to parabolic induction from two Speh representations}
	
	\author{David Ginzburg}
	
	\address{School of Mathematical Sciences, Sackler Faculty of Exact Sciences, Tel-Aviv University, Israel
		69978} \email{ginzburg@tauex.tau.ac.il}

	\thanks{This research was supported by the ISRAEL SCIENCE FOUNDATION
		(grant No. 722/26).}
	
	\author{David Soudry}
	\address{School of Mathematical Sciences, Sackler Faculty of Exact Sciences, Tel-Aviv University, Israel
		69978} \email{soudry@tauex.tau.ac.il}
	
	
	
	
	
	\maketitle
	
	\begin{abstract}
		We prove identities of Siegel-Weil type for Eisenstein series on general linear groups, attached to parabolic induction from two Speh representations corresponding to the same irreducible, cuspidal representation. Such an identity relates a residue, at a (simple) pole of an Eisenstein series, as above, to a value of a similarly structured Eisenstein series at a point of holomorphy.
	\end{abstract}

	\section{Introduction}
	
	Let $F$ be a number field and $\BA$ its ring of adeles. Let $\tau$ be an irreducible, unitary, cuspidal, automorphic representation of $\GL_n(\BA)$. For a positive integer $k$, denote by $\Delta(\tau,k)$ the corresponding Speh representation of $\GL_{kn}(\BA)$. Let $m_1, m_2$ be two positive integers. Denote $m=m_1+m_2$ and $\underline{m}=(m_1,m_2)$. Let $\underline{s}=(s_1,s_2)\in \BC^2$. Consider the parabolic induction
    \begin{equation}\label{1}
      \rho_{\Delta(\tau,\underline{m}),\underline{s}}=\Delta(\tau,m_1)|\cdot|^{s_1}\times \Delta(\tau,m_2)|\cdot|^{s_2}.
      \end{equation}
  Note that 
  \begin{equation}\label{2}
  \rho_{\Delta(\tau,\underline{m}),\underline{s}}=|\det\cdot|^{\frac{s_1+s_2}{2}}\Delta(\tau,m_1)|\cdot|^{\frac{s}{2}}\times \Delta(\tau,m_2)|\cdot|^{-\frac{s}{2}},
  \end{equation}
  	where $s=s_1-s_2$. Consider an Eisenstein series attached to $\rho_{\Delta(\tau,\underline{m}),\underline{s}}$. Due to \eqref{2}, we may consider it as a meromorphic function of the single variable $s$. Then we know (see \cite{GS25}),
  \begin{thm}\label{thm A}
  	The list of poles of the Eisenstein series attached to $\rho_{\Delta(\tau,\underline{m}),\underline{s}}$, for $Re(s)\geq 0$, is
  	$$
  	s=\frac{m}{2}-i, \ \ 0\leq i\leq min(m_1,m_2)-1.
  	$$
  	All these poles are simple.
  \end{thm}
Denote
\begin{equation}\label{3}
	\rho_{\Delta(\tau,\underline{m}),s}=\Delta(\tau,m_1)|\cdot|^{\frac{s}{2}}\times \Delta(\tau,m_2)|\cdot|^{-\frac{s}{2}},
\end{equation}
and for a smooth, holomorphic section $f_{\Delta(\tau,\underline{m}),s}$ of $\rho_{\Delta(\tau,\underline{m}),s}$, denote the corresponding Eisenstein series by $E(f_{\Delta(\tau,\underline{m}),s})$.	
The main result of this paper is the following identity, of Siegel-Weil type, showing that each residue in the last theorem is a value of a similar Eisenstein series at a point of holomorphy.
\begin{thm}\label{thm B}
	Let $f_{\Delta(\tau,\underline{m}),s}$ be a smooth, holomorphic section of $\rho_{\Delta(\tau,\underline{m}),s}$.\\
	1. Assume that $m_1\geq m_2$. Let $0\leq i\leq m_2-1$. Then there is a smooth, holomorphic section $\varphi_{\Delta(\tau,(m-i,i)),\zeta}$ of $\rho_{\Delta(\tau,(m-i,i)),\zeta}$, such that
	\begin{equation}\label{4}
		Res_{s=\frac{m}{2}-i}E(f_{\Delta(\tau,\underline{m}),s})=value_{\zeta=\frac{m_1-m_2}{2}}E(\varphi_{\Delta(\tau,(m-i,i)),\zeta}).
	\end{equation}
2. Assume that $m_1<m_2$. Let $0\leq i\leq m_1-1$. Then there is a smooth, holomorphic section $\varphi_{\Delta(\tau,(i,m-i)),\zeta}$ of $\rho_{\Delta(\tau,(i,m-i)),\zeta}$, such that
\begin{equation}\label{5}
	Res_{s=\frac{m}{2}-i}E(f_{\Delta(\tau,\underline{m}),s})=value_{\zeta=\frac{m_2-m_1}{2}}E(\varphi_{\Delta(\tau,(i,m-i)),\zeta}).
\end{equation}
\end{thm}
Referring to the identities \eqref{4}, \eqref{5}, as of Siegel-Weil type, is after the work of Kudla and Rallis \cite{KR94}. One of their main theorems (\cite{KR94}, p. 5) relates leading terms of certain degenerate Eisenstein series on a symplectic group, at certain poles. In detail, consider an Eisenstein series on $\Sp_{2n}(\BA)$, attached to the parabolic induction $|\det_{\GL_r}\cdot|^s\rtimes 1_{\Sp_{2(n-r)}}$. Denote it by $E_r^{2n}(f_s)$, 
for a smooth holomorphic section $f_s$. Assume that $\frac{n+1}{2}<r\leq n+1$. Then $E_r^{2n}(f_s)$ has a double pole at $s=\frac{r-1}{2}$, as the section varies. Denote the leading term of the corresponding Laurent expansion by $R^{(-2)}_{s=\frac{r-1}{2}}(E_r^{2n}(f_s))$. Consider the Eisenstein series $E_{n-r+1}^{2n}(f'_z)$. Then it has a simple pole at $z=\frac{n-r}{2}$, as the section varies. Denote the corresponding residue by $R^{(-1)}_{z=\frac{n-r}{2}}(E_{n-r+1}^{2n}(f'_z))$. Then, for compatible sections,
\begin{equation}\label{6}
R^{(-2)}_{s=\frac{r-1}{2}}(E_r^{2n}(f_s))=R^{(-1)}_{z=\frac{n-r}{2}}(E_{n-r+1}^{2n}(f'_z)).
\end{equation}
Consider the Eisenstein series $E_n^{2n}(f''_\zeta)$, the degenerate Siegel Eisenstein series. Then it has a simple pole at $\zeta=r-\frac{n+1}{2}$, as the section varies. Kudla and Rallis prove also that, for compatible sections,
\begin{equation}\label{7}
R^{(-1)}_{z=\frac{n-r}{2}}(E_{n-r+1}^{2n}(f'_z))=R^{(-1)}_{\zeta=r-\frac{n+1}{2}}(E_n^{2n}(f''_\zeta)).
\end{equation}
The relations in Theorem \ref{thm B} bear a similarity with \eqref{6}. They compare leading terms of two Eisenstein series of similar structure at two points, where the order of the l.h.s. is one less than that of the r.h.s.. In \eqref{6}, a double pole is related to a simple pole, and in Theorem \ref{thm B}, a simple pole is related to a value. We note that the work of Kudla, Rallis was generalized to other classical groups, by Ikeda, Moeglin, Ichino, Yamana, Gan-Qiu-Takeda and others. See [GQT14] and the references therein.

We would like to present another point of view, namely that of the theorem of Franke, \cite{F98}. It says that an automorphic form is a linear combination of values of Eisenstein series and their derivatives, where the Eisenstein series are attached to parabolic inductions from square-integrable automorphic representations. Theorem \ref{thm B} shows that derivatives do not contribute to the r.h.s. of \eqref{4}, \eqref{5}, and, of course, it identifies the value. As for $R^{(-2)}_{s=\frac{r-1}{2}}(E_r^{2n}(f_s))$, it is equal by \eqref{7} to $R^{(-1)}_{\zeta=r-\frac{n+1}{2}}(E_n^{2n}(f''_\zeta))$, which is square-integrable (\cite{KR94}, p. 8). Again, this fits Franke's theorem, and with no contribution from derivatives.

\section{The constant term of $E(f_{\Delta(\tau,\underline{m}),\underline{s}})$ along $V_{k,mn-k}$}

Let $0<k<mn$ be an integer. We compute the constant term
\begin{equation}\label{1.1}
	E^{V_{k,mn-k}}(f_{\Delta(\tau,\underline{m}),\underline{s}})(g)=\int_{[V_{k,mn-k}]}E(f_{\Delta(\tau,\underline{m}),\underline{s}})(ug)du.
\end{equation}
Here, $V_{\underline{m}n}$ is the unipotent radical of the standard parabolic subgroup $P_{\underline{m}n}$ of $\GL_{mn}$, with Levi part, $M_{\underline{m}n}$, isomorphic to $\GL_{m_1n}\times \GL_{m_2n}$. We denote $[V_{k,mn-k}]=V_{k,mn-k}(F)\backslash V_{k,mn-k}(\BA)$.

\begin{prop}\label{prop 1.1}
	The constant term \eqref{1.1} is zero, unless $k$ is a multiple of $n$. Let $k=\ell n$, $1\leq \ell<m$. Then
\begin{equation}\label{1.1'}
E^{V_{\ell n,(m-\ell)n}}(f_{\Delta(\tau,\underline{m}),\underline{s}})(g)=\sum_j\sum_{\gamma_1,\gamma_2}M_{w_{jn},\underline{s}}(f^{V_{(j,m_1-j)n}\times V_{(\ell-j,m_2-\ell+j)n}}_{\Delta(\tau,\underline{m}),\underline{s}})(\begin{pmatrix}\gamma_1\\&\gamma_2\end{pmatrix}g),
\end{equation}
where the first sum runs over non-negative integers $j$, such that $\ell-m_2\leq j\leq \ell, m_1$; the second sum is over $\gamma_1\in P_{(j,\ell-j)n}(F)\backslash \GL_{\ell n}(F)$,\\
$\gamma_2\in P_{(m_1-j)n,(m_2-\ell+j)n}(F)\backslash \GL_{(m-\ell)n}(F)$; $M_{w_{jn},\underline{s}}$ is the intertwining operator corresponding to the Weyl element $w_{jn}$ written below in \eqref{1.3}. See \eqref{1.8} below. Finally, 	
$$
f^{V_{(j,m_1-j)n}\times V_{(\ell-j,m_2-\ell+j)n}}_{\Delta(\tau,\underline{m}),\underline{s}}(g)=\int f_{\Delta(\tau,\underline{m}),\underline{s}}(\begin{pmatrix}v_1\\&v_2\end{pmatrix}g)d(v_1,v_2),
$$
where the integration is over $[V_{(j,m_1-j)n}\times V_{(\ell-j,m_2-\ell+j)n}]$.
This is the composition of the section with the application of the constant terms to $\Delta(\tau,m_1)$, along $V_{(j,m_1-j)n}$ and to $\Delta(\tau,m_2)$, along $V_{(\ell-j,m_2-\ell+j)n}$. 
\end{prop}

\begin{proof}

	The computation is standard. Assume that $Re(s)$ is large, so that $E(f_{\Delta(\tau,\underline{m}),\underline{s}})$ is given by the absolutely convergent series
$$
E(f_{\Delta(\tau,\underline{m}),\underline{s}})(g)=\sum_{\gamma\in P_{\underline{m}n}(F)\backslash \GL_{mn}(F)}f_{\Delta(\tau,\underline{m}),\underline{s}}(\gamma g).
$$
Then the constant term \eqref{1.1} is equal to
\begin{equation}\label{1.2}
\sum_{w\in P_{\underline{m}n}(F)\backslash \GL_{mn}(F)/P_{k,mn-k}(F)}	\int_{[V_{k,mn-k}]}\sum_{\eta\in Q^{(w)}(F)\backslash P_{k,mn-k}(F)}f_{\Delta(\tau,\underline{m}),\underline{s}})(w\eta ug)du,
\end{equation}
where
$$
Q^{(w)}=P_{k,mn-k}\cap w^{-1}P_{\underline{m}n}w.
$$
Here is a set of representatives for the double cosets $P_{\underline{m}n}\backslash \GL_{mn}/P_{k,mn-k}$,
\begin{equation}\label{1.3}
	w_r=\begin{pmatrix}I_r&0&0&0\\0&0&I_{m_1n-r}&0\\0&I_{k-r}&0&0\\0&0&0&I_{m_2n-k+r}\end{pmatrix},\ k-m_2n\leq r\leq k, m_1n.
\end{equation}    	
The subgroup $Q^{(w_r)}$ consists of the elements
\begin{equation}\label{1.4}
	\begin{pmatrix} a_1&a_2&x_1&x_2\\0&a_4&0&x_4\\&&b_1&b_2\\&&0&b_4\end{pmatrix},\ a_1\in \GL_r,\ a_4\in \GL_{k-r}, b_1\in \GL_{m_1n-r}.
\end{equation}	
The conjugation $g\mapsto w_rgw_r^{-1}$ takes the element \eqref{1.4} to
\begin{equation}\label{1.5}
	\begin{pmatrix}a_1&x_1&a_2&x_2\\0&b_1&0&b_2\\&&a_4&x_4\\&&0&b_4\end{pmatrix}.
\end{equation}
Continue unfolding \eqref{1.2} by factoring $P_{k,mn-k}(F)$ from the right by $V_{k,mn-k}(F)$. Fix $r$ as above, and take the following representatives of  $Q^{(w_r)}(F)\backslash P_{k,mn-k}(F)/V_{k,mn-k}(F)$,
$$
\begin{pmatrix}\gamma_1\\&\gamma_2\end{pmatrix},\ \gamma_1\in P_{r,k-r}(F)\backslash \GL_k(F),\ \gamma_2\in P_{m_1n-r,m_2n-k+r}(F)\backslash \GL_{mn-k}(F).
$$  	
The integral in \eqref{1.2}, corresponding to $w_r$ is then equal to
\begin{equation}\label{1.6}
	\sum	
    	\int_{V_{k,mn-k}^{(r)}(F)\backslash V_{k,mn-k}(\BA)}f_{\Delta(\tau,\underline{m}),\underline{s}}(w_rv\begin{pmatrix}\gamma_1\\&\gamma_2\end{pmatrix}g)dv,
    \end{equation}
where the summation is over 
$$
\gamma_1\in P_{r,k-r}(F)\backslash \GL_k(F), \gamma_2\in P_{m_1n-r,m_2n-k+r}(F)\backslash \GL_{mn-k}(F),
$$
and
$$
V_{k,mn-k}^{(r)}=V_{k,mn-k}\cap Q^{(w_r)}.
$$
Note that by \eqref{1.5}, the subgroup $V_{k,mn-k}^{(r)}$ conjugates by $w_r$ to the subgroup of elements   	
\begin{equation}\label{1.7}
	\begin{pmatrix}I_r&x_1&0&x_2\\0&I_{m_1n-r}&0&0\\&&I_{k-r}&x_4\\&&0&I_{m_2n-k+r}\end{pmatrix}.
\end{equation}    	
Thus, when we factor the integration in the integral in \eqref{1.6} through $V_{k,mn-k}^{(r)}(\BA)$, we get 
\begin{equation}\label{1.8}
	\int_{M_{(k-r)\times (m_1n-r)}(\BA)}f^{V_{r,m_1n-r}\times V_{k-r,m_2n-k+r}}_{\Delta(\tau,\underline{m}),\underline{s}}(w_ru_{r,k}(x)\begin{pmatrix}\gamma_1\\&\gamma_2\end{pmatrix}g)dv,
\end{equation}    	
 where 
 $$
 u_{r,k}(x)=\begin{pmatrix}I_r\\&I_{k-r}&x\\&&I_{m_1n-r}\\&&&&I_{m_2n-k+r}\end{pmatrix},
 $$
 and
 $$
 f^{V_{r,m_1n-r}\times V_{k-r,m_2n-k+r}}_{\Delta(\tau,\underline{m}),\underline{s}}(g)=\int_{ [V_{r,m_1n-r}\times V_{k-r,m_2n-k+r}]}f_{\Delta(\tau,\underline{m}),\underline{s}}(\begin{pmatrix}v_1\\&v_2\end{pmatrix}g)d(v_1,v_2).
 $$
 Thus, we need to apply to $\Delta(\tau,m_1)$ the constant term along $V_{r,m_1n-r}$, and to $\Delta(\tau,m_2)$ the constant term along $V_{k-r,m_2n-k+r}$. Since $\tau$ is cuspidal, this constant term is zero, unless $k$ and $r$ are multiples of $n$. So, assume that $k=\ell n$, $r=jn$, where $1\leq \ell< m$, $\ell-m_2\leq j\leq \ell, m_1$. Denote the integral \eqref{1.8} by
  	$M_{w_{jn},\underline{s}}(f^{V_{(j,m_1-j)n}\times V_{(\ell-j,m_2-\ell+j)n}}_{\Delta(\tau,\underline{m}),\underline{s}})(g)$. This is the intertwining operator corresponding to $w_{jn}$. Now \eqref{1.1'} follows from \eqref{1.2}, \eqref{1.6}, \eqref{1.8}.  We proved \eqref{1.1'} for $Re(s)$ sufficiently large. By analytic continuation \eqref{1.1'} holds as an equality of meromorphic functions.
\end{proof}    	
    	
Denote by $\Delta(\tau,m_1)^{V_{(j,m_1-j)n}}$ the space of constant terms along $V_{(j,m_1-j)n}$, applied to the elements of $\Delta(\tau,m_1)$. When we restrict these constant terms to the Levi subgroup $M_{(j,m_1-j)n}(\BA)$, we get 
\begin{equation}\label{1.9}
\Delta(\tau,m_1)^{V_{(j,m_1-j)n}}_{|M_{(j,m_1-j)n}(\BA)}=\delta^{\frac{1}{2}}_{P_{(j,m_1-j)n}}\Delta(\tau,j)|\cdot|^{-\frac{m_1-j}{2}}\otimes \Delta(\tau,m_1-j)|\cdot|^{\frac{j}{2}}.
\end{equation}
As functions on $\GL_{m_1n}(\BA)$,
\begin{equation}\label{1.10}
\Delta(\tau,m_1)^{V_{(j,m_1-j)n}}\subset \Delta(\tau,j)|\cdot|^{-\frac{m_1-j}{2}}\times \Delta(\tau,m_1-j)|\cdot|^{\frac{j}{2}}.
\end{equation}
The r.h.s. denotes the (normalized) parabolic induction.	Similarly,
\begin{equation}\label{1.11}
	\Delta(\tau,m_2)^{V_{(\ell-j,m_2-\ell+j)n}}_{|M_{(\ell-j,m_2-\ell+j)n}(\BA)}=\delta^{\frac{1}{2}}_{P_{(\ell-j,m_2-\ell+j)n}}\Delta(\tau,\ell-j)|\cdot|^{-\frac{m_2-\ell+j}{2}}\otimes \Delta(\tau,m_2-\ell+j)|\cdot|^{\frac{\ell-j}{2}}.
\end{equation}
As functions on $\GL_{m_2n}(\BA)$,
\begin{equation}\label{1.12}
	\Delta(\tau,m_2)^{V_{(\ell-j,m_2-\ell+j)n}}\subset \Delta(\tau,\ell-j)|\cdot|^{-\frac{m_2-\ell+j}{2}}\times \Delta(\tau,m_2-\ell+j)|\cdot|^{\frac{\ell-j}{2}}.
\end{equation}	    	
Thus, by \eqref{1.9}-\eqref{1.12}, as a function on $\GL_{mn}(\BA)$,\\
\\
$ f^{V_{(j,m_1-j)n}\times V_{(\ell-j,m_2-\ell+j)n}}_{\Delta(\tau,\underline{m}),\underline{s}}\in$
\\
\begin{equation}\label{1.13}
 \Delta(\tau,j)|\cdot|^{s_1-\frac{m_1-j}{2}}\times \Delta(\tau,m_1-j)|\cdot|^{s_1+\frac{j}{2}}\times \Delta(\tau,\ell-j)|\cdot|^{s_2-\frac{m_2-\ell+j}{2}}\times \Delta(\tau,m_2-\ell+j)|\cdot|^{s_2+\frac{\ell-j}{2}}.
\end{equation}
We view $M_{w_{jn},\underline{s}}$ as an operator on the parabolic induction \eqref{1.13}. In fact, it is the intertwining operator on \eqref{1.13}, which flips the two middle factors. This is clear from \eqref{1.5}, \eqref{1.8}. Thus,\\
\\
$ M_{w_{jn},\underline{s}}(f^{V_{(j,m_1-j)n}\times V_{(\ell-j,m_2-\ell+j)n}}_{\Delta(\tau,\underline{m}),\underline{s}})\in$
\\
\begin{equation}\label{1.14}
	\Delta(\tau,j)|\cdot|^{s_1-\frac{m_1-j}{2}}\times \Delta(\tau,\ell-j)|\cdot|^{s_2-\frac{m_2-\ell+j}{2}}\times \Delta(\tau,m_1-j)|\cdot|^{s_1+\frac{j}{2}}\times \Delta(\tau,m_2-\ell+j)|\cdot|^{s_2+\frac{\ell-j}{2}}.
\end{equation}	    	
The poles (in $s=s_1-s_2$) of $M_{w_{jn},\underline{s}}$, as an operator on \eqref{1.13} are the same as those of the intertwining operator
\begin{equation}\label{1.15}
	M^{m_1-j,\ell-j}_{(s_1+\frac{j}{2},s_2-\frac{m_2-\ell+j}{2})}: \Delta(\tau,m_1-j)|\cdot|^{s_1+\frac{j}{2}}\times \Delta(\tau,\ell-j)|\cdot|^{s_2-\frac{m_2-\ell+j}{2}}\rightarrow 
\end{equation}
$$
	 \Delta(\tau,\ell-j)|\cdot|^{s_2-\frac{m_2-\ell+j}{2}}\times \Delta(\tau,m_1-j)|\cdot|^{s_1+\frac{j}{2}},
$$
flipping the two factors. At a place $v$ of $F$, let $M^{m_1-j,\ell-j}_{(s_1+\frac{j}{2},s_2-\frac{m_2-\ell+j}{2});v}$ be the similar local intertwining operator at $v$, and let $N^{m_1-j,\ell-j}_{(s_1+\frac{j}{2},s_2-\frac{m_2-\ell+j}{2});v}$ be the corresponding normalized local intertwining operator.
\begin{lem}\label{lem 1.2}
For each place $v$, the normalized local intertwining operator\\ 	$N^{m_1-j,\ell-j}_{(s_1+\frac{j}{2},s_2-\frac{m_2-\ell+j}{2});v}$ is holomorphic and nonzero when $Re(s)\geq 0$. In this domain, the poles of $	M^{m_1-j,\ell-j}_{(s_1+\frac{j}{2},s_2-\frac{m_2-\ell+j}{2})}$ are the same as the poles of
\begin{equation}\label{1.16}
	\prod_{i=o}^{\ell-j-1}\frac{L(\tau\times\hat{\tau},s+j+\frac{m_2-m_1}{2}-i)}{L(\tau\times\hat{\tau},s+\frac{m}{2}-i)}.
\end{equation}
These are the poles of $ M_{w_{jn},\underline{s}}$, as an operator on \eqref{1.13}, when $Re(s)\geq 0$.
\end{lem}
\begin{proof}
The statement about the local intertwining operator follows from the fact that $$
Re((s_1+\frac{j}{2})-(s_2-\frac{m_2-\ell+j}{2}))\geq 0.
$$
Indeed, the l.h.s. is equal to $Re(s)+\frac{j}{2}+\frac{j-(\ell-m_2)}{2}$. This is non-negative, since $Re(s),j\geq 0$, and $j\geq \ell-m_2$. Now, the statement follows from \cite{MW89}, I.10. We have\\
\\
$	M^{m_1-j,\ell-j}_{(s_1+\frac{j}{2},s_2-\frac{m_2-\ell+j}{2})}=$

\begin{equation}\label{1.17}
 \frac{L(\Delta(\tau,m_1-j)\times \Delta(\hat{\tau},\ell-j),s+\frac{j}{2}+\frac{m_2-\ell+j}{2})}{\epsilon L(\Delta(\tau,m_1-j)\times \Delta(\hat{\tau},\ell-j),s+1+\frac{j}{2}+\frac{m_2-\ell+j}{2})}\prod_vN^{m_1-j,\ell-j}_{(s_1+\frac{j}{2},s_2-\frac{m_2-\ell+j}{2});v}.
\end{equation}
Here, $\epsilon=\epsilon(\Delta(\tau,m_1-j)\times \Delta(\hat{\tau},\ell-j),s+\frac{j}{2}+\frac{m_2-\ell+j}{2},\psi)$, the epsilon factor (with respect to a nontrivial character $\psi$ of $F\backslash \BA$). This is an exponential and hence holomorphic and nonzero everywhere.
Note that the product over $v$ in \eqref{1.17}, when applied to a given section, becomes a finite product. A simple calculation shows that the ratio of $L$-functions in \eqref{1.17} is equal to \eqref{1.16}.
\end{proof} 	
By \eqref{1.14}, we can write, for $a\in \GL_{\ell n}(\BA)$, $b\in \GL_{(m-\ell)n}(\BA)$,\\
\\
$ M_{w_{jn},\underline{s}}(f^{V_{(j,m_1-j)n}\times V_{(\ell-j,m_2-\ell+j)n}}_{\Delta(\tau,\underline{m}),\underline{s}})(\begin{pmatrix}a\\&b\end{pmatrix})=$
\\
\begin{equation}\label{1.18}
\sum \xi_{\Delta(\tau,j)|\cdot|^{s_1-\frac{m_1-j}{2}}\times \Delta(\tau,\ell-j)|\cdot|^{s_2-\frac{m_2-\ell+j}{2}}}(a)\varphi_{\Delta(\tau,m_1-j)|\cdot|^{s_1+\frac{j}{2}}\times \Delta(\tau,m_2-\ell+j)|\cdot|^{s_2+\frac{\ell-j}{2}}}(b),
\end{equation}
where $\xi_{\Delta(\tau,j)|\cdot|^{s_1-\frac{m_1-j}{2}}\times \Delta(\tau,\ell-j)|\cdot|^{s_2-\frac{m_2-\ell+j}{2}}}$, $\varphi_{\Delta(\tau,m_1-j)|\cdot|^{s_1+\frac{j}{2}}\times \Delta(\tau,m_2-\ell+j)|\cdot|^{s_2+\frac{\ell-j}{2}}}$ are meromorphic sections of the indicated parabolic inductions.	This, Proposition \ref{prop 1.1}, \eqref{1.6} and \eqref{1.8}, show 
\begin{prop}\label{prop 1.3}
Let $1\leq \ell<m$. Then, for $a\in \GL_{\ell n}(\BA)$, $b\in \GL_{(m-\ell)n}(\BA)$,\\
\\
$	E^{V_{(\ell,m-\ell)n}}(f_{\Delta(\tau,\underline{m}),\underline{s}})(\begin{pmatrix}a\\&b\end{pmatrix})=\sum_j\sum$
\begin{equation}\label{1.19}
	E(\xi_{\Delta(\tau,j)|\cdot|^{s_1-\frac{m_1-j}{2}}\times \Delta(\tau,\ell-j)|\cdot|^{s_2-\frac{m_2-\ell+j}{2}}})(a)E(\varphi_{\Delta(\tau,m_1-j)|\cdot|^{s_1+\frac{j}{2}}\times \Delta(\tau,m_2-\ell+j)|\cdot|^{s_2+\frac{\ell-j}{2}}})(b),
\end{equation}
where $j$ runs over non-negative integers $j$, such that $\ell-m_2\leq j\leq \ell, m_1$. For each $j$, the second sum runs over the meromorphic sections as in \eqref{1.18}.
\end{prop}
Let us examine the exponents of our Eisenstein series on $GL_{mn}(\BA)$ along $P_{(\ell,m-\ell)n}$, that is we consider the action of the center of the Levi subgroup $M_{(\ell,m-\ell)n}(\BA)$ on each summand in \eqref{1.19} corresponding to $j$ in the first sum.
\begin{lem}\label{lem 1.4}
	Let $j_1, j_2$ be two different, non-negative integers, such that $\ell-m_2\leq j_1, j_2\leq \ell, m_1$. Then their corresponding exponents (in \eqref{1.19}) are equal iff
	\begin{equation}\label{1.20}
		j_1+j_2=\frac{m_1-m_2}{2}+\ell-s.
	\end{equation}
\end{lem}
\begin{proof}
Each term, corresponding to an index $j$, in \eqref{1.19}, has the central character
\begin{equation}\label{1.21}	
	\omega_\tau(\lambda\mu)|\lambda|^{nj(s_1-\frac{m_1-j}{2})+n(\ell-j)(s_2-\frac{m_2-\ell+j}{2})}|\mu|^{n(m_1-j)(s_1+\frac{j}{2})+n(m_2-\ell+j)(s_2+\frac{\ell-j}{2})},
\end{equation}
where $\omega_\tau$ is the central character of $\tau$. Thus, the exponents, for $j_1\neq j_2$, are equal iff
$$
(j_1(s_1-\frac{m_1-j_1}{2})+(\ell-j_1)(s_2-\frac{m_2-\ell+j_1}{2}),(m_1-j_1)(s_1+\frac{j_1}{2})+(m_2-\ell+j_1)(s_2+\frac{\ell-j_1}{2}))=
$$
$$
(j_2(s_1-\frac{m_1-j_2}{2})+(\ell-j_2)(s_2-\frac{m_2-\ell+j_2}{2}),(m_1-j_2)(s_1+\frac{j_2}{2})+(m_2-\ell+j_2)(s_2+\frac{\ell-j_2}{2})).
$$
Compare the first coordinates. We get
$$
(j_1-j_2)s+(j_1-j_2)\frac{m_2-m_1}{2}+(j_1-j_2)\frac{j_1+j_2}{2}-(j_1-j_2)\frac{2\ell-(j_1+j_2)}{2}=0.
$$
Since $j_1\neq j_2$, this is the same as 
$$
s+\frac{m_2-m_1}{2}+j_1+j_2-\ell=0,
$$
as we want. When we compare the second coordinates, we get the same condition.
\end{proof}
\begin{cor}\label{cor 1.5}
	For each exponent corresponding to an index $j$ in \eqref{1.19}, there is at most one other index $j'$ with the same exponent.
\end{cor}
\section{Proof of Theorem \ref{thm B} when $m_1\geq m_2$}

Let $m_1\geq m_2$ and $1\leq \ell \leq m_2$. Consider the constant term \eqref{1.1}, for $k=\ell n$, 
 $	E^{V_{(\ell,m-\ell)n}}(f_{\Delta(\tau,\underline{m}),\underline{s}})$. Consider in \eqref{1.13} the two indices $j=\ell,0$. For $j=\ell$, we view the function on
  $\GL_{(m-\ell)n}(\BA)$, $g\mapsto f^{V_{(\ell,m_1-\ell)n}\times I_{m_2n}}_{\Delta(\tau,\underline{m}),\underline{s}}(\begin{pmatrix}I_{\ell n}\\&g\end{pmatrix})$, as an element of
   $\Delta(\tau,\ell)|\cdot|^{s_1-\frac{m_1-\ell}{2}}\times  (\Delta(\tau,m_1-\ell)|\cdot|^{s_1+\frac{\ell}{2}}\times \Delta(\tau,m_2)|\cdot|^{s_2})$, evaluated at 
   $\begin{pmatrix}I_{\ell n}\\&g\end{pmatrix}$. Let us further evaluate this element at $(I_{\ell n}, I_{(m-\ell)n})$. Denote this by
\begin{equation}\label{2.1}
    f^{V_{(\ell,m_1-\ell)n}\times I_{m_2n}}_{\Delta(\tau,\underline{m}),\underline{s}}(\begin{pmatrix}I_{\ell n}\\&g\end{pmatrix})(I_{\ell n},I_{(m-\ell)n})
    =|\det g|^{-\frac{\ell n}{2}}f'_{\Delta(\tau,m_1-\ell)|\cdot|^{s_1+\frac{\ell}{2}}\times \Delta(\tau,m_2)|\cdot|^{s_2}}(g),  
     \end{equation}
     where $f'_{\Delta(\tau,m_1-\ell)|\cdot|^{s_1+\frac{\ell}{2}}\times \Delta(\tau,m_2)|\cdot|^{s_2}}$ is a smooth, holomorphic section of\\
      $\Delta(\tau,m_1-\ell)|\cdot|^{s_1+\frac{\ell}{2}}\times \Delta(\tau,m_2)|\cdot|^{s_2}$. Note that $M_{w_{\ell n,\underline{s}}}=id$.
     Let $j=0$. In \eqref{1.14}, we view the function on $\GL_{(m-\ell)n}(\BA)$,
     $g\mapsto M_{w_{0,\underline{s}}}(f^{I_{m_1 n},V_{(\ell,m_2-\ell)n}}_{\Delta(\tau,\underline{m}),\underline{s}}(\begin{pmatrix}I_{\ell n}\\&g\end{pmatrix})$, as an element of $\Delta(\tau,\ell)|\cdot|^{s_2-\frac{m_2-\ell}{2}}\times  (\Delta(\tau,m_1)|\cdot|^{s_1}\times \Delta(\tau,m_2-\ell)|\cdot|^{s_2+\frac{\ell}{2}})$, evaluated at 
     	$\begin{pmatrix}I_{\ell n}\\&g\end{pmatrix}$. Again, evaluate this element at $(I_{\ell n}, I_{(m-\ell)n})$. Denote this by
     	\begin{equation}\label{2.2}
     	M_{w_{0,\underline{s}}}(f^{I_{m_1n}\times V_{(\ell,m_2-\ell)n}}_{\Delta(\tau,\underline{m})),\underline{s}}(\begin{pmatrix}I_{\ell n}\\&g\end{pmatrix})(I_{\ell n,(m-\ell)n})
     	=
     	\end{equation}
     	$$
     	|\det g|^{-\frac{\ell n}{2}}f''_{\Delta(\tau,m_1)|\cdot|^{s_1}}\times \Delta(\tau,m_2-\ell)|\cdot|^{s_2+\frac{\ell}{2}}(g),  
     	$$
     	where $f''_{\Delta(\tau,m_1)|\cdot|^{s_1}\times \Delta(\tau,m_2-\ell)|\cdot|^{s_2+\frac{\ell}{2}}}$	is a smooth,
     		 meromorphic section of $\Delta(\tau,m_1)|\cdot|^{s_1}\times \Delta(\tau,m_2-\ell)|\cdot|^{s_2+\frac{\ell}{2}}$. 
     
     \begin{prop}\label{prop 2.1}
     Let $m_1> m_2$ and $1\leq \ell \leq m_2$. Assume that $\ell\leq m_1-m_2$. Let $\underline{s}$ be such that $s=\frac{m_1-m_2}{2}$. Then, in the notation above, as the section $f_{\Delta(\tau,\underline{m}),\underline{s}}$ varies, $E(f'_{\Delta(\tau,m_1-\ell)|\cdot|^{s_1+\frac{\ell}{2}}\times \Delta(\tau,m_2)|\cdot|^{s_2}})$ has a simple pole at $s=\frac{m_1-m_2}{2}$, as well as $f''_{\Delta(\tau,m_1)|\cdot|^{s_1}\times \Delta(\tau,m_2-\ell)|\cdot|^{s_2+\frac{\ell}{2}}}$, and we have
     \begin{equation}\label{2.3}
     Res_{s=\frac{m_1-m_2}{2}}E(f'_{\Delta(\tau,m_1-\ell)|\cdot|^{s_1+\frac{\ell}{2}}\times \Delta(\tau,m_2)|\cdot|^{s_2}})=
     	\end{equation}
     	$$
     	-value_{s=\frac{m_1-m_2}{2}}E(Res_{s=\frac{m_1-m_2}{2}}f''_{\Delta(\tau,m_1)|\cdot|^{s_1}\times \Delta(\tau,m_2-\ell)|\cdot|^{s_2+\frac{\ell}{2}}}).
     	$$
     	\end{prop}
     	\begin{proof}
     	We first note that $E(f_{\Delta(\tau,\underline{m}),\underline{s}})$ is holomorphic at $s=\frac{m_1-m_2}{2}$. This follows from Theorem \ref{thm A}. We have
     	$$
     	s_1+\frac{\ell}{2}-s_2=\frac{m_1-m_2+\ell}{2}=\frac{m-\ell}{2}-(m_2-\ell).
     	$$
     	Since $m_1<m_2$, $m_2-\ell<m_1-\ell$. Of course, $m_2-\ell<m_2$. By Theorem \ref{thm A}, $E(f'_{\Delta(\tau,m_1-\ell)|\cdot|^{s_1+\frac{\ell}{2}}\times \Delta(\tau,m_2)|\cdot|^{s_2}})$ has a simple pole at $s=\frac{m_1-m_2}{2}$, as the section $f_{\Delta(\tau,\underline{m}),\underline{s}}$ varies. We have 
     	$$
     	s_1-(s_2-\frac{m_2-\ell}{2})=\frac{m_1-\ell}{2}> 0.
     	$$
     	By Lemma \ref{lem 1.2}, for $j=0$, we get that $f''_{\Delta(\tau,m_1)|\cdot|^{s_1}\times \Delta(\tau,m_2-\ell)|\cdot|^{s_2+\frac{\ell}{2}}}$ has a simple pole at $s=\frac{m_1-m_2}{2}$. 
     	We have
     	$$
     	s_1-(s_2+\frac{\ell}{2})=\frac{m_1-m_2-\ell}{2}\geq 0,
     	$$
     	by our assumption. Since
     	$$
     	s_1-(s_2+\frac{\ell}{2})=\frac{m-\ell}{2}-m_2,
     	$$
     	$E(\varphi_{\Delta(\tau,m_1)|\cdot|^{s_1}\times \Delta(\tau,m_2-\ell)|\cdot|^{s_2+\frac{\ell}{2}}})$, is holomorphic at $s=\frac{m_1-m_2}{2}$, for any smooth section $\varphi_{\Delta(\tau,m_1)|\cdot|^{s_1}\times \Delta(\tau,m_2-\ell)|\cdot|^{s_2+\frac{\ell}{2}}}$, which is holomorphic at $s=\frac{m_1-m_2}{2}$. Thus, $value_{s=\frac{m_1-m_2}{2}}E(Res_{s=\frac{m_1-m_2}{2}}f''_{\Delta(\tau,m_1)|\cdot|^{s_1}\times \Delta(\tau,m_2-\ell)|\cdot|^{s_2+\frac{\ell}{2}}})$ makes sense. Now consider \eqref{1.19} in Proposition \ref{prop 1.3}. The l.h.s. is holomorphic at $s=\frac{m_1-m_2}{2}$, but there are terms in the r.h.s. of \eqref{1.19} which have a pole at this point. In particular, the sum of all these terms, with the same exponent, along the center of $M_{(\ell,m-\ell)}(\BA)$, is holomorphic at this point. Now, the proposition follows from Lemma \ref{lem 1.4}, with $j_1=\ell$, $j_2=0$, $s=\frac{m_1-m_2}{2}$.

     	\end{proof}
     
 		{\bf Proof of Theorem \ref{thm B}, when $m_1\geq m_2$:} 
 		Assume that $m_1\geq m_2$. Substitute $m_1+\ell$ instead of $m_1$ in Proposition \ref{prop 2.1}, and let $1\leq \ell\leq m_2$. Note that $m_1+\ell>m_2$, and $\ell\leq (m_1+\ell)-m_2$. Next, substitute $s_1-\frac{\ell}{2}$ instead of $s_1$ in Proposition \ref{prop 2.1}. The l.h.s. of \eqref{2.3} becomes
 		\begin{equation}\label{2.4}
 		Res_{s_1-s_2=\frac{m_1+m_2}{2}-(m_2-\ell)}E(f'_{\Delta(\tau,m_1)|\cdot|^{s_1}\times \Delta(\tau,m_2)|\cdot|^{s_2}}).
 		\end{equation}
 		The r.h.s. becomes
 		\begin{equation}\label{2.5}
 		-value_{s_1-s_2=\frac{m_1+m_2}{2}-(m_2-\ell)}E(\xi_{\Delta(\tau,m_1+\ell)|\cdot|^{s_1-\frac{\ell}{2}}\times \Delta(\tau, m_2-\ell)|\cdot|^{s_2+\frac{\ell}{2}}}),
 		\end{equation}
 		where $\xi_{\Delta(\tau,m_1+\ell)|\cdot|^{s_1-\frac{\ell}{2}}\times \Delta(\tau, m_2-\ell)|\cdot|^{s_2+\frac{\ell}{2}}}$ is any holomorphic section, such that its value at $\underline{s}$, with $s_1-s_2=\frac{m_1+m_2}{2}-(m_2-\ell)$ is equal to 
 		$$
 		Res_{s=\frac{m_1+m_2}{2}-(m_2-\ell)}f''_{\Delta(\tau,m_1+\ell)|\cdot|^{s_1-\frac{\ell}{2}}\times \Delta(\tau,m_2-\ell)|\cdot|^{s_2+\frac{\ell}{2}}}.
 		$$
 		denote $i=m_2-\ell$. Then $0\leq i\leq m_2-1$. Then
 		\begin{equation}\label{2.6}
 		Res_{s_1-s_2=\frac{m}{2}-i}E(f'_{\Delta(\tau,m_1)|\cdot|^{s_1}\times \Delta(\tau,m_2)|\cdot|^{s_2}})=
 		\end{equation}
 		$$
 			-value_{s_1-s_2=\frac{m}{2}-i}E(\xi_{\Delta(\tau,m-i)|\cdot|^{s_1-\frac{m_2-i}{2}}\times \Delta(\tau, i)|\cdot|^{s_2+\frac{m_2-i}{2}}}).
 		$$
 		Using \eqref{2}, rewrite \eqref{2.6} as
 		\begin{equation}\label{2.7}
 		|\det\cdot|^{\frac{s_1+s_2}{2}}\cdot Res_{s=\frac{m}{2}-i}E(f'_{\Delta(\tau,m_1)|\cdot|^{\frac{s}{2}}\times \Delta(\tau,m_2)|\cdot|^{-\frac{s}{2}}})=
 		\end{equation}
 		$$
 		-|\det\cdot|^{\frac{s_1+s_2}{2}}\cdot value_{s=\frac{m}{2}-i}E(\xi_{\Delta(\tau,m-i)|\cdot|^{\frac{s-m_2+i}{2}}\times \Delta(\tau, i)|\cdot|^{-\frac{s-m_2+i}{2}}}).
 		$$
 		Hence
 		\begin{equation}\label{2.8}
 		 Res_{s=\frac{m}{2}-i}E(f'_{\Delta(\tau,m_1)|\cdot|^{\frac{s}{2}}\times \Delta(\tau,m_2)|\cdot|^{-\frac{s}{2}}})=
 			-value_{\zeta=\frac{m_1-m_2}{2}}E(\xi_{\Delta(\tau,m-i)|\cdot|^{\frac{\zeta}{2}}\times \Delta(\tau, i)|\cdot|^{-\frac{\zeta}{2}}}).
 		\end{equation}
 		Finally, let us recall how the sections in \eqref{2.8} are obtained. We assumed that $m_1\geq m_2$, and $0\leq i\leq m_2-1$. Let $f_{\Delta(\tau,(m-i,m_2)),(z_1,z_2)}$ be a smooth, holomorphic section of $\Delta(\tau,m-i)|\cdot|^{z_1}\times \Delta(\tau,m_2)|\cdot|^{z_2}$. Then, for $g\in \GL_{mn}(\BA)$, 
 		\begin{equation}\label{2.9}
 		f'_{\Delta(\tau,m_1)|\cdot|^{s_1}\times \Delta(\tau,m_2)|\cdot|^{s_2}}(g)=
 		\end{equation}
 		$$
 		|\det g|^{\frac{(m_2-i) n}{2}}f^{V_{(m_2-i,m_1)n\times I_{m_2n}}}_{\Delta(\tau,(m-i,m_2)),(s_1-\frac{m_2-i}{2},s_2)}(\begin{pmatrix}I_{(m_2-i)n}\\&g\end{pmatrix})(I_{(m_2-i)n},I_{mn}).
 		$$
 		See \eqref{1.13}, \eqref{2.1}. Clearly, as the section $f_{\Delta(\tau,(m-i,m_2)),(z_1,z_2)}$ varies,\\ $f'_{\Delta(\tau,m_1)|\cdot|^{s_1}\times \Delta(\tau,m_2)|\cdot|^{s_2}}$ describes an arbitrary smooth, holomorphic section of\\ $\Delta(\tau,m_1)|\cdot|^{s_1}\times \Delta(\tau,m_2)|\cdot|^{s_2}$. 
 		
 		The section $\xi_{\Delta(\tau,m-i)|\cdot|^{z_1}\times \Delta(\tau,i)|\cdot|^{z_2}}$ is any holomorphic section, such that
 		\begin{equation}\label{2.10}
 		value_{s=\frac{m}{2}-i}\xi_{\Delta(\tau,m-i)|\cdot|^{s_1-\frac{m_2-i}{2}}\times \Delta(\tau,i)|\cdot|^{s_2+\frac{m_2-i}{2}}}(g)=|\det g|^{\frac{(m_2-i) n}{2}}
 		\end{equation}
 		$$
 		Res_{s=\frac{m}{2}-i}M_{w_{0,(s_1-\frac{m_2-i}{2},s_2)}}(f^{I_{m_1n}\times V_{(m_2-i,i)n}}_{\Delta(\tau,(m-i,m_2)),(s_1-\frac{m_2-i}{2},s_2)})(\begin{pmatrix}I_{(m_2-i)n}\\&g\end{pmatrix})(I_{(m_2-i)n,mn}).
 		$$
 		See \eqref{1.14}, \eqref{2.2}, \eqref{2.5}. Note that the representation $$\Delta(\tau,m-i)|\cdot|^{\frac{m_1-m_2}{4}}\times \Delta(\tau,i)|\cdot|^{-\frac{m_1-m_2}{4}}$$ is irreducible. 
 		Indeed, the segments $[\frac{1-m+i}{2}+\frac{m_1-m_2}{4},\frac{m-i-1}{2}+\frac{m_1-m_2}{4}]$ and\\ $[\frac{1-i}{2}-\frac{m_1-m_2}{4},\frac{i-1}{2}-\frac{m_1-m_2}{4}]$ are not linked. Theorem \ref{thm B} is proved for $m_1\geq m_2$.
 			 
\section{Proof of Theorem \ref{thm B} when $m_1<m_2$}
	
	Assume that $m_1<m_2$, and let $2m_1\leq \ell<m$. Assume also that $\ell\geq m_2$. Consider in \eqref{1.13} the two indices $j=m_1, \ell-m_2$. Note that, for $j=m_1$, $M_{w_{m_1n,\underline{s}}}=id$. We view the function on $\GL_{\ell n}(\BA)$, $g\mapsto f^{I_{m_1n}\times V_{(\ell-m_1,m-\ell)n}}_{\Delta(\tau,\underline{m}),\underline{s}}(\begin{pmatrix}g\\&I_{(m-\ell)n}\end{pmatrix})$, as an element of $\Delta(\tau,m_1)|\cdot|^{s_1}\times \Delta(\tau,\ell-m_1)|\cdot|^{s_2-\frac{m-\ell}{2}}\times \Delta(\tau,m-\ell)|\cdot|^{s_2+\frac{\ell-m_1}{2}}$ evaluated at $\begin{pmatrix}g\\&I_{(m-\ell)n}\end{pmatrix}$. As in \eqref{2.1}, let us further evaluate at $(I_{\ell n}, I_{(m-\ell)n})$. Denote this by
	\begin{equation}\label{3.1}
	 f^{I_{m_1n}\times V_{(\ell-m_1,m-\ell)n}}_{\Delta(\tau,\underline{m}),\underline{s}}(\begin{pmatrix}g\\&I_{(m-\ell)n}\end{pmatrix})(I_{\ell n}, I_{(m-\ell)n})=
	 \end{equation}
	 $$
	 |det g|^{\frac{(m-\ell)n}{2}}f'_{\Delta(\tau,m_1)|\cdot|^{s_1}\times \Delta(\ell-m_1)|\cdot|^{s_2-\frac{m-\ell}{2}}}(g),
	 $$
	 where  $f'_{\Delta(\tau,m_1)|\cdot|^{s_1}\times \Delta(\ell-m_1)|\cdot|^{s_2-\frac{m-\ell}{2}}}$ is a smooth, holomorphic section of\\
	  $\Delta(\tau,m_1)|\cdot|^{s_1}\times \Delta(\tau,\ell-m_1)|\cdot|^{s_2-\frac{m-\ell}{2}}$. Let $j=\ell-m_2$. We view the function on $\GL_{\ell n}(\BA)$, $g\mapsto M_{w_{(\ell-m_2)n},\underline{s}}(f_{\Delta(\tau,\underline{m}),\underline{s}}^{V_{(\ell-m_2,m-\ell)n}\times I_{m_2n}})(\begin{pmatrix}g\\&I_{(m-\ell)n}\end{pmatrix})$, as an element of $\Delta(\tau, \ell-m_2)|\cdot|^{s_1-\frac{m-\ell}{2}}\times \Delta(\tau,m_2)|\cdot|^{s_2}\times \Delta(\tau,m-\ell)|\cdot|^{s_1+\frac{\ell-m_2}{2}}$ evaluated at $\begin{pmatrix}g\\&I_{(m-\ell)n}\end{pmatrix}$. Evaluate this at $(I_{\ell n},I_{(m-\ell)n})$. Denote this by
	 \begin{equation}\label{3.2}
	 M_{w_{(\ell-m_2)n},\underline{s}}(f_{\Delta(\tau,\underline{m}),\underline{s}}^{V_{(\ell-m_2,m-\ell)}\times I_{m_2n}})(\begin{pmatrix}g\\&I_{(m-\ell)n}\end{pmatrix})(I_{\ell n},I_{(m-\ell)n})=
	 \end{equation}
	 $$
	 |\det g|^{\frac{(m-\ell)n}{2}}f''_{\Delta(\tau, \ell-m_2)|\cdot|^{s_1-\frac{m-\ell}{2}}\times \Delta(\tau,m_2)|\cdot|^{s_2}}(g),
	 $$
	 where $f''_{\Delta(\tau, \ell-m_2)|\cdot|^{s_1-\frac{m-\ell}{2}}\times \Delta(\tau,m_2)|\cdot|^{s_2}}$ is a smooth, meromorphic section of $\Delta(\tau, \ell-m_2)|\cdot|^{s_1-\frac{m-\ell}{2}}\times \Delta(\tau,m_2)|\cdot|^{s_2}$.	
    	
    \begin{prop}\label{prop 3.1}
    	Let $m_1< m_2$ and $2m_1\leq \ell < m$. Assume that $\ell\geq m_2$. Let $\underline{s}$ be such that $s=\frac{m_2-m_1}{2}$. Then, in the notation above, as the section $f_{\Delta(\tau,\underline{m}),\underline{s}}$ varies, $E(f'_{\Delta(\tau,m_1)|\cdot|^{s_1}\times \Delta(\tau,\ell-m_1)|\cdot|^{s_2-\frac{m-\ell}{2}}})$ has a simple pole at $s=\frac{m_2-m_1}{2}$, as well as $f''_{\Delta(\tau,\ell-m_2)|\cdot|^{s_1-\frac{m-\ell}{2}}\times \Delta(\tau,m_2)|\cdot|^{s_2}}$, and we have
    	\begin{equation}\label{3.3}
    		Res_{s=\frac{m_2-m_1}{2}}E(f'_{\Delta(\tau,m_1)|\cdot|^{s_1}\times \Delta(\tau,\ell-m_1)|\cdot|^{s_2-\frac{m-\ell}{2}}})=
    	\end{equation}
    	$$
    	-value_{s=\frac{m_2-m_1}{2}}E(Res_{s=\frac{m_2-m_1}{2}}f''_{\Delta(\tau,\ell-m_2)|\cdot|^{s_1-\frac{m-\ell}{2}}\times \Delta(\tau,m_2)|\cdot|^{s_2}}).
    	$$
    \end{prop}

\begin{proof}
	As is Proposition \ref{prop 2.1}, we first note that $E(f_{\Delta(\tau,\underline{m}),\underline{s}})$ is holomorphic at $s=\frac{m_2-m_1}{2}$. This follows from Theorem \ref{thm A}. We have
	$$
	s_1-(s_2-\frac{m-\ell}{2})=\frac{m_2-m_1}{2}+\frac{m_1+m_2-\ell}{2}=m_2-\frac{\ell}{2}=\frac{\ell}{2}-(\ell-m_2).
	$$
	By our assumptions, $0\leq \ell-m_2<m_1, \ell-m_1$. By Theorem \ref{thm A},\\
	 $E(f'_{\Delta(\tau,m_1)|\cdot|^{s_1}\times \Delta(\tau,\ell-m_1)|\cdot|^{s_2-\frac{m-\ell}{2}}})$ has a simple pole at $s=\frac{m_2-m_1}{2}$, as the section $f_{\Delta(\tau,\underline{m}),\underline{s}}$ varies. We have
	 $$
	 s_1+\frac{\ell-m_2}{2}-s_2=s+\frac{\ell-m_2}{2}=\frac{\ell-m_1}{2}>0.
	 $$
	 Next, in Lemma \ref{lem 1.2}, for $j=\ell-m_2$, the product \eqref{1.16} is
	 $$
	 \frac{\prod_{i=0}^{m_2-1}L(\tau\times\hat{\tau},s+\ell-\frac{m}{2}-i)}{\prod_{i=0}^{m_2-1}L(\tau\times\hat{\tau},s+\frac{m}{2}-i)}.
	 	$$
	 	The numerator has a double pole at $s=\frac{m_2-m_1}{2}$, and the denominator has a simple pole. All in all, this has a simple pole at $s=\frac{m_2-m_1}{2}$. By Lemma \ref{lem 1.2}, $f''_{\Delta(\tau,\ell-m_2)|\cdot|^{s_1-\frac{m-\ell}{2}}\times \Delta(\tau,m_2)|\cdot|^{s_2}}$ has a simple pole at $s=\frac{m_2-m_1}{2}$. We have
	 	$$
	 	s_1-\frac{m-\ell}{2}-s_2=\frac{m_2-m_1}{2}-\frac{m-\ell}{2}=\frac{\ell}{2}-m_1\geq 0,
	 	$$
	 	by our assumption. Since $m_1>\ell-m_2$, $E(\varphi_{\Delta(\tau,\ell-m_2)|\cdot|^{s_1-\frac{m-\ell}{2}}\times \Delta(\tau,m_2)|\cdot|^{s_2}})$ is holomorphic at $s=\frac{m_2-m_1}{2}$, for any smooth section $\varphi_{\Delta(\tau,\ell-m_2)|\cdot|^{s_1-\frac{m-\ell}{2}}\times \Delta(\tau,m_2)|\cdot|^{s_2}}$, which is holomorphic at $s=\frac{m_2-m_1}{2}$, so that\\ $value_{s=\frac{m_2-m_1}{2}}E(Res_{s=\frac{m_2-m_1}{2}}f''_{\Delta(\tau,\ell-m_2)|\cdot|^{s_1-\frac{m-\ell}{2}}\times \Delta(\tau,m_2)|\cdot|^{s_2}})$ is meaningful. 
	 	As in Proposition \ref{prop 2.1}, consider \eqref{1.19} in Proposition \ref{prop 1.3}. The l.h.s. is holomorphic at $s=\frac{m_2-m_1}{2}$, but there are terms in the r.h.s. of \eqref{1.19} which have a pole at this point. In particular, the sum of all these terms, with the same exponent, along the center of $M_{(\ell,m-\ell)}(\BA)$, is holomorphic at this point. Now, the proposition follows from Lemma \ref{lem 1.4}, with $j_1=m_1$, $j_2=\ell-m_2$, $s=\frac{m_2-m_1}{2}$.

\end{proof}

    {\bf Proof of Theorem \ref{thm B}, when $m_2> m_1$:} Assume that $m_2> m_1$. Let $1\leq k\leq m_1$ be an integer. Substitute $m_2+k$ instead of $m_2$ in Proposition \ref{prop 3.1}. Of course, $m_1<m_2+k$. Let $\ell=m_1+m_2=m$. Note that $2m_1<\ell<m+k$, and $\ell\geq m_2+k$. Then
    \begin{equation}\label{3.4}
    	Res_{s=\frac{m_2+k-m_1}{2}}E(f'_{\Delta(\tau,m_1)|\cdot|^{s_1}\times \Delta(\tau,m_2)|\cdot|^{s_2-\frac{k}{2}}})=
    \end{equation}
    $$
    -value_{s=\frac{m_2+k-m_1}{2}}E(\xi_{\Delta(\tau,m_1-k)|\cdot|^{s_1-\frac{k}{2}}\times \Delta(\tau,m_2+k)|\cdot|^{s_2}}),
    $$
    where $\xi_{\Delta(\tau,m_1-k)|\cdot|^{s_1-\frac{k}{2}}\times \Delta(\tau,m_2+k)|\cdot|^{s_2}}$ is any holomorphic section, such that its value at $\underline{s}$, with $s_1-s_2=\frac{m_2+k-m_1}{2}$, is equal to 
    $$
    Res_{s=\frac{m_2+k-m_1}{2}}f''_{\Delta(\tau,m_1-k)|\cdot|^{s_1-\frac{k}{2}}\times \Delta(\tau,m_2+k)|\cdot|^{s_2}}.
    $$
    Now, substitute in \eqref{3.4} $s_2+\frac{k}{2}$ instead of $s_2$, and denote $i=m_1-k$. Then $0\leq i\leq m_1-1$, and
\begin{equation}\label{3.5}
	Res_{s=\frac{m}{2}-i}E(f'_{\Delta(\tau,m_1)|\cdot|^{s_1}\times \Delta(\tau,m_2)|\cdot|^{s_2}})=
\end{equation}
$$
-value_{s=\frac{m}{2}-i}E(\xi_{\Delta(\tau,i)|\cdot|^{s_1-\frac{m_1-i}{2}}\times \Delta(\tau,m-i)|\cdot|^{s_2+\frac{m_1-i}{2}}}).
$$
As in \eqref{2.7}, \eqref{2.8}, we get
\begin{equation}\label{3.6}
Res_{s=\frac{m}{2}-i}E(f'_{\Delta(\tau,m_1)|\cdot|^{\frac{s}{2}}\times \Delta(\tau,m_2)|\cdot|^{-\frac{s}{2}}})=-value_{\zeta=\frac{m_2-m_1}{2}}E(\xi_{\Delta(\tau,i)|\cdot|^{\frac{\zeta}{2}}\times \Delta(\tau,m-i)|\cdot|^{-\frac{\zeta}{2}}}).
	\end{equation}
Finally, let us recall how the sections in \eqref{3.6} are obtained. Let $m_1<m_2$ and $0\leq i\leq m_1-1$.	Let $f_{\Delta(\tau,(m_1,m-i)),(z_1,z_2)}$ be a smooth, holomorphic section of $\Delta(\tau,m_1)|\cdot|^{z_1}\times \Delta(\tau,m-i)|\cdot|^{z_2}$. Then, for $g\in \GL_{mn}(\BA)$, 
\begin{equation}\label{3.7}
	f'_{\Delta(\tau,m_1)|\cdot|^{s_1}\times \Delta(\tau,m_2)|\cdot|^{s_2}}(g)=
\end{equation}
$$
|\det g|^{-\frac{(m_1-i) n}{2}}f^{I_{m_1n}\times V_{(m_2,m_1-i)n}}_{\Delta(\tau,(m_1,m-i)),(s_1,s_2+\frac{m_1-i}{2})}(\begin{pmatrix}g\\&I_{(m_1-i)n}\end{pmatrix})(I_{mn},I_{(m_1-i)n}).
$$
See \eqref{1.13}, \eqref{3.1}. Clearly, as the section $f_{\Delta(\tau,(m_1,m-i)),(z_1,z_2)}$ varies,\\ $f'_{\Delta(\tau,m_1)|\cdot|^{s_1}\times \Delta(\tau,m_2)|\cdot|^{s_2}}$ describes an arbitrary smooth, holomorphic section of\\ $\Delta(\tau,m_1)|\cdot|^{s_1}\times \Delta(\tau,m_2)|\cdot|^{s_2}$.
	
The section $\xi_{\Delta(\tau,i)|\cdot|^{z_1}\times \Delta(\tau,m-i)|\cdot|^{z_2}}$ is any holomorphic section, such that
\begin{equation}\label{3.8}
	value_{s=\frac{m}{2}-i}\xi_{\Delta(\tau,i)|\cdot|^{s_1-\frac{m_1-i}{2}}\times \Delta(\tau,m-i)|\cdot|^{s_2+\frac{m_1-i}{2}}}(g)=|\det g|^{-\frac{(m_1-i) n}{2}}
\end{equation}
$$
Res_{s=\frac{m}{2}-i}M_{w_{in,(s_1,s_2+\frac{m_1-i}{2})}}(f_{\Delta(\tau,(m_1,m-i)),(s_1,s_2+\frac{m_1-i}{2})}^{V_{(i,m_1-i)n}\times I_{m_2n}})(\begin{pmatrix}g\\&I_{(m_i-i)n}\end{pmatrix})(I_{mn}, I_{(m_1-i)n}).
$$
See \eqref{1.14}, \eqref{3.2}, \eqref{3.5}. Note that the representation $$\Delta(\tau,i)|\cdot|^{\frac{m_2-m_1}{4}}\times \Delta(\tau,m-i)|\cdot|^{-\frac{m_2-m_1}{4}}$$ is irreducible, since $[\frac{1-i}{2}+\frac{m_2-m_1}{4},\frac{i-1}{2}+\frac{m_2-m_1}{4}]$, $[\frac{1-m+i}{2}-\frac{m_2-m_1}{4},\frac{m-i-1}{2}-\frac{m_2-m_1}{4}]$ are not linked. Theorem \ref{thm B} is proved for $m_1< m_2$.

  \end{document}